\documentclass [12pt,a4paper] {amsart}

\usepackage [T1] {fontenc}
\usepackage {amssymb}
\usepackage {amsmath}
\usepackage {amsthm}
\usepackage{tikz}
\usepackage{hyperref}
\usepackage{url}
\usetikzlibrary{arrows,decorations.pathmorphing,backgrounds,positioning,fit,petri,cd}
\usepackage{caption}
\usepackage{subcaption}

 \usepackage[titletoc,toc,title]{appendix}
 
\usepackage[margin=0.8in]{geometry}

\usetikzlibrary{arrows,decorations.pathmorphing,backgrounds,positioning,fit,petri}

\usetikzlibrary{matrix}

\usepackage{enumitem}

\DeclareMathOperator {\tp} {tp}

\DeclareMathOperator {\acl} {acl}
\DeclareMathOperator {\dcl} {dcl}

\DeclareMathOperator {\op} {op}

\DeclareMathOperator {\ECBDA} {ECBDA}
\DeclareMathOperator {\ECDA} {ECDA}

\theoremstyle {definition}
\newtheorem {definition}{Definition} [section]
\newtheorem{example} [definition] {Example}

\newtheorem* {notation} {Notation}

\newtheorem* {observation} {Observation}

\theoremstyle {plain}

\newtheorem {lemma} [definition] {Lemma}
\newtheorem {theorem} [definition] {Theorem}
\newtheorem {proposition} [definition] {Proposition}
\newtheorem {corollary} [definition] {Corollary}

\theoremstyle {remark}
\newtheorem {remark} [definition] {Remark}

\makeatletter
\newcommand {\forksym} {\raise0.2ex\hbox{\ooalign{\hidewidth$\vert$\hidewidth\cr\raise-0.9ex\hbox{$\smile$}}}}
\def\@forksym@#1#2{\mathrel{\mathop{\forksym}\displaylimits_{#2}}}
\def\forkind{\@ifnextchar_{\@forksym@}{\forksym}}

\makeatother

\begin {document}

\title{Existentially closed De Morgan algebras}

\author{Vahagn Aslanyan}
\address{Department of Mathematical Sciences, CMU, Pittsburgh, PA 15213, USA}
\address{Institute of Mathematics of the National Academy of Sciences, Yerevan 0019, Armenia}

\address{\texttt{Current address} School of Mathematics, University of East Anglia, Norwich, NR4 7TJ, UK}
\email{V.Aslanyan@uea.ac.uk}



\maketitle

\begin{abstract}
We show that the theory of De Morgan algebras has a model completion and axiomatise it. Then we prove that it is $\aleph_0$-categorical and describe definable and algebraic closures in that theory. We also obtain similar results for Boole-De Morgan algebras.
\end{abstract}


\section{Introduction}

A \emph{De Morgan algebra} \cite{moisil,birkhoff,Gratzer-lattice-theory,kalman,balbes-dwinger} is a structure ${D}:=(D;+,\cdot, ~\bar{}~, 0, 1)$ where $(D;+,\cdot, 0,1)$ is a bounded distributive lattice with a largest element $1$ and a smallest element $0$ and $~\bar{}~$ is a unary operation satisfying the following identities:
$$\forall x (\bar{\bar{x}}=x) \mbox{ and }  \forall x, y (\overline{x+y} = \bar{x}\cdot \bar{y}).$$ If, in addition, $\forall x (x+\bar{x}=1)$ then ${D}$ is said to be a \emph{Boolean algebra}. Note that henceforth we will denote the unary operation of a Boolean algebra (Boolean negation) by $'$ and that of a De Morgan algebra (De Morgan negation) by $~\bar{}~$. It is customary to use $\vee$ and $\wedge$ for lattice operations (known as \textit{join} and \textit{meet}), but we use $+$ and $\cdot$ instead to avoid confusion with logical disjunction and conjunction which are used throughout the paper. 

De Morgan algebras are important structures in algebra and mathematical logic and have various applications. In particular, they are closely related to algebraic logic and, more specifically, to Belnap--Dunn logic and Relevance logic \cite{Ross-Belnap,Ross-belnap-dunn}, and fuzzy logic \cite{hajek-fuzzy}. Note that the standard fuzzy algebra $([0,1]; \min(x,y), \max(x,y), 1-x, 0, 1)$ is a De Morgan algebra. Furthermore, De Morgan algebras have applications to multi-valued simulations of digital circuits \cite{Brzozowski-hazard}.

It is well known that a Boolean algebra is existentially closed if and only if it is atomless. The theory consisting of the axioms of Boolean algebras and the additional axiom 
\begin{equation*}
    \forall x (x>0 \rightarrow \exists y (0<y<x))
\end{equation*}
is the model completion of the theory of Boolean algebras and has quantifier elimination \cite{Poizat}. Similarly, a bounded distributive lattice is existentially closed if and only if it is atomless and complemented \cite{Schmid-ec-dl}. The latter means that every element in the lattice has a complement, that is, $$\forall x \exists y (xy=0 \wedge x+y=1).$$ Complements in distributive lattices are unique, hence a bounded complemented distributive lattice is just the underlying lattice of a Boolean algebra.

In this paper we give a first-order characterisation of existentially closed De Morgan algebras and thus obtain a model companion of the theory of De Morgan algebras. We also observe that De Morgan algebras satisfy the amalgamation property, hence the aforementioned theory is actually the model completion of the theory of De Morgan algebras and has quantifier elimination. Unlike the case of Boolean algebras, the theory that we obtain for existentially closed De Morgan algebras is somewhat more complicated. In particular we have an axiom scheme stating that certain systems of one-variable equations and inequations have solutions. 

Further, we study model theoretic properties of that theory. In particular, we show that it is $\aleph_0$-categorical, i.e. it has a unique countable model up to isomorphism (which is also the case for atomless Boolean algebras), and describe definable and algebraic closures in existentially closed De Morgan algebras. Actually we prove all those results for Boole-De Morgan algebras first (those are lattices equipped with a Boolean negation and a De Morgan negation, see Section \ref{Boole-DeMorgan-section}) and then translate them to the language of De Morgan algebras. This is possible due to the fact that the underlying lattice of an existentially closed De Morgan algebra is complemented.

Note that a Boolean algebra is uniquely determined by its lattice structure. Unlike this, we have a lot of freedom in defining a De Morgan structure on a given distributive lattice. 
It is this fact that makes the model theoretic treatment of De Morgan algebras significantly harder than that of Boolean algebras.

After finishing the work on this paper it was brought to my attention by James Raftery that \cite{Clark-Davey} discusses similar questions and may have some overlap with the current paper. Indeed, \cite[Chapter 5]{Clark-Davey} discusses existentially and algebraically closed algebras in varieties from the point of view of natural dualities (Theorem 5.3.5). Then the authors give the aforementioned characterisation of existentially and algebraically closed Boolean algebras and bounded distributive lattices (Theorems 5.4.1 and 5.4.2). Further, they ask in Exercise 5.10 to show that a De Morgan algebra is existentially closed if and only if it is atomless and complemented. However, we claim that this description is incorrect, and we will prove this in the next section (Remark \ref{wrong-proof}). Actually, if it were true then every atomless Boolean algebra would be existentially closed in the variety of De Morgan algebras, which is not the case. Moreover, it is easy to see that one has to impose some conditions on the negation (unary operation) of a De Morgan algebra to make it existentially closed; properties of the underlying lattice alone cannot characterise existentially closed De Morgan algebras. Our characterisation is somewhat more complicated and we believe there is no simple description as in the case of Boolean algebras. One might also get a characterisation of existentially closed De Morgan algebras using the methods of \cite{Clark-Davey} but apparently one has to do a substantial amount of work for that. On the other hand, our proof presented here is quite elementary and does not require any advanced theory.

\section{Preliminaries}

\subsection{Model theoretic preliminaries}

In this section we recall some basic model theoretic notions that will be used throughout the paper. More notions and results will be recalled later in the paper when we need them. The reader is referred to \cite{Mar,TZ} for details.

\begin{definition}
Given two structures $A \subseteq B$ (in the same language), $A$ is called \emph{existentially closed} in $B$ if every existential formula with parameters from $A$ that is true in $B$ is also true in $A$ (in other words, if a quantifier-free formula with parameters from $A$ has a realisation in $B$ then it also has a realisation in $A$). A model of a theory $T$ is said to be \emph{existentially closed} if it is existentially closed in all extensions which are also models of $T$.
\end{definition}

\begin{definition}
Let $T$ be a first-order theory.

\begin{itemize}
    \item $T$ is \emph{model complete} if all models of $T$ are existentially closed. Equivalently, $T$ is model complete if every formula is equivalent to an existential formula modulo $T$.
    \item A theory $T'$ is a \emph{model companion} of $T$ if $T'$ is model complete and every model of $T$ can be embedded into a model of $T'$ and vice versa. 
    \item $T$ has the \emph{amalgamation property} if for any models $M_0, M_1, M_2$ of $T$ with embeddings $f_i: M_0 \hookrightarrow M_i,~ i=1,2,$ there is a model $M$ of $T$ with embeddings $g_i:M_i \hookrightarrow M$ such that $g_1 \circ f_1 = g_2 \circ f_2$.
    \item If $T$ has the amalgamation property then a model companion of $T$ is also called a \emph{model completion}.
\end{itemize}
\end{definition}

\begin{remark}
A theory has at most one model companion (hence at most one model completion). So we may speak of \emph{the} model companion (completion) of a theory assuming it exists.
\end{remark}

\begin{remark}
The model completion of a universal theory (if it exists) admits elimination of quantifiers, that is, every formula is equivalent to a quantifier-free formula modulo that theory (the model completion). 
\end{remark}

Let us give a few examples. 
\begin{example}
\begin{itemize}
    \item The theory of algebraically closed fields is the model completion of the theory of fields.
    \item  The theory of real closed fields is the model completion of the theory of ordered fields (in the language of ordered rings).
    \item The theory of atomless Boolean algebras is the model completion of the theory of Boolean algebras.
\end{itemize}
\end{example}

\subsection{Distributive lattices and Boolean algebras}

Given a bounded distributive lattice ${L}=(L;+,\cdot, 0,1)$, its \emph{dual} lattice is the lattice given by the reverse order on $L$. It is denoted by ${L}^{\op}:= (L;\cdot,+, 1,0)$ where the sequence of operations suggests that the meet of ${L}$ is the join of ${L}^{\op}$ and vice versa. The operations on the direct product ${L} \times {L}^{\op}:= (L \times L; +, \cdot, 0, 1)$ are defined by
$$(x_1,x_2)+ (y_1,y_2) = (x_1+y_1, x_2\cdot y_2),~ (x_1,x_2)\cdot (y_1,y_2) = (x_1\cdot y_1, x_2 + y_2),~ 0 = (0,1),~ 1 = (1,0).$$
It can be made into a De Morgan algebra by defining $\overline{(x,y)} = (y,x)$. Furthermore, every De Morgan algebra can be embedded into such a one. Given an arbitrary De Morgan algebra ${D}=(D;+,\cdot, ~\bar{}~, 0, 1)$, the map $i:x\mapsto (x,\bar{x})$ is a De Morgan embedding of ${D}$ into ${D}_L \times {D}_L^{\op}$ where ${D}_L$ is the underlying lattice of ${D}$. Note that this construction is sometimes referred to as \emph{twist-product} (see, for example, \cite{twist-product}) and is attributed to Kalman \cite{kalman}.


The following results will be used in our analysis of existentially closed De Morgan algebras.

\begin{theorem}[{\cite[Chapter 6]{Poizat}}]
The theory of atomless Boolean algebras has quantifier elimination and is the model completion of the theory of Boolean algebras.
\end{theorem}

\begin{theorem}[\cite{Schmid-ec-dl}]
The theory of bounded, atomless and complemented distributive lattices is the model companion of the theory of bounded distributive lattices.
\end{theorem}

In fact, the theory of (bounded) distributive lattices has the amalgamation property (see \cite{Gratzer-lattice-theory}) and since it is a universal theory, we get the following consequence.
\begin{corollary}
The theory of bounded, atomless and complemented distributive lattices is the model completion of the theory of bounded distributive lattices. It is complete and admits quantifier elimination.
\end{corollary}

\begin{remark}\label{wrong-proof}
Now we prove that it is not true that a De Morgan algebra with a complemented and atomless underlying lattice is existentially closed. Indeed, this would imply that an atomless Boolean algebra is existentially closed as a De Morgan algebra. Pick such an algebra ${B}=(B;+,\cdot, ~\bar{}~, 0, 1)$ and embed it into ${D}:={B}_L \times {B}_L^{\op}$ as above (where ${B}_L$ is the underlying lattice of ${B}$). Note that ${D}$ is not a Boolean algebra and we claim that ${B}$ is not existentially (even algebraically) closed in it. To this end observe that the sentence $\exists x (\bar{x} = x)$ is true in ${D}$ (e.g. for any element $u \in B$ we have $\overline{(u,u)} = (u,u)$ in ${D}$) but not in ${B}$.
\end{remark}


\section{Boole-De Morgan algebras}\label{Boole-DeMorgan-section}

In this section we introduce Boole-De Morgan algebras and make a few observations about them which will be used later. We refer the reader to \cite{Mov-Asl-Boole-de} for more details.

\begin{definition}
A \emph{Boole-De Morgan algebra} is an algebra $(A;+,\cdot, ', ~\bar{}~, 0, 1)$ where $(A;+,\cdot, ',  0, 1)$ is a Boolean algebra and $(A;+,\cdot, ~\bar{}~, 0, 1)$ is a De Morgan algebra.
\end{definition}

\begin{observation}
The Boolean and De Morgan negations commute in a Boole-De Morgan algebra, that is, $\overline{(x')} = (\bar{x})'$. To prove this notice that 
$$\overline{(x')} + \bar{x} = \overline{x'\cdot x} = 1,~ \overline{(x')} \cdot \bar{x} = \overline{x'+x} = 0.$$
\end{observation}

If ${B}=(B;+, \cdot, ', 0,1)$ is a Boolean algebra then the direct product ${B}_L \times {B}_L^{\op}$ is a Boole-De Morgan algebra where the Boolean and De Morgan negations are defined as follows:
$$(x,y)' = (x',y'),~ \overline{(x,y)} = (y,x).$$ Moreover, for every Boole-De Morgan algebra ${B}$ the map $x \mapsto (x,\bar{x})$ is an embedding ${B} \hookrightarrow {B}_L \times {B}_L^{\op}$ due to the above observation.

\begin{figure}

  \centering
\begin{tikzpicture}[auto, node distance=2cm, every loop/.style={},
                    thick,main node/.style={font=\sffamily\large}]

  \node[main node] (1) {$1$};
  \node[main node] (2) [below left of=1] {$a$};
  \node[main node] (3) [below right of=2] {$0$};
  \node[main node] (4) [below right of=1] {$b$};

  \path[every node/.style={font=\sffamily\small}]
    (1) edge node [left] {} (4)
      
    (2) edge node [right] {} (1)

    (3) edge node [right] {} (2)
        
    (4) edge node [left] {} (3);
\end{tikzpicture}
  \caption{The underlying lattice of $\mathbf{4}$}
  \label{fig1}
\end{figure}
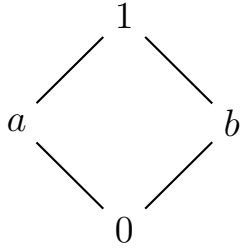

The lattice given in Figure \ref{fig1} can be made into a Boole-De Morgan algebra by defining $a'=b,~ b'=a, ~ \bar{a}=a,~ \bar{b}=b,~ 0'=\bar{0}=1,~ 1' = \bar{1} = 0$. It will be denoted by $\mathbf{4}:=(\{ 0,1, a, b\}; +, \cdot, ', ~\bar{}~, 0, 1)$. Its subalgebra with domain $\{ 0, 1\}$ is denoted by $\mathbf{2}$.

It is proved in \cite[Theorem 3.12]{Mov-Asl-Boole-de} that $\mathbf{2}$ and $\mathbf{4}$ are subdirectly irreducible and in fact these are the only subdirectly irreducible Boole-De Morgan algebras (this is an analogue of the well known characterisation of subdirectly irreducible De Morgan algebras \cite{kalman}). In particular, we get the following embedding theorem.

\begin{proposition}[\cite{Mov-Asl-Boole-de}]\label{BDM-embedding-4}
Every Boole-De Morgan algebra can be embedded into a direct power of $\mathbf{4}$.
\end{proposition}

\begin{corollary}\label{local-finite}
Boole-De Morgan algebras are locally finite, i.e. every finitely generated Boole-De Morgan algebra is finite.
\end{corollary}

\begin{lemma}
Let ${D} = (D;+,\cdot, ~\bar{}~, 0, 1)$ be an existentially closed De Morgan algebra. Then $(D;+,\cdot, 0, 1)$ is atomless and complemented.
\end{lemma}
\begin{proof}
Extend $(D;+,\cdot, 0, 1)$ to an atomless and complemented distributive lattice $(L;+, \cdot, 0, 1)$. The map $x\mapsto (x, \bar{x})$ gives an embedding of the De Morgan algebra ${D}$ into $L \times L^{\op}$. The latter is a De Morgan algebra with atomless and complemented underlying lattice. Thus, ${D}$ can be embedded into an atomless and complemented De Morgan algebra and since it is existentially closed, it must itself be atomless and complemented.
\end{proof}

\begin{corollary}
Let $(D;+,\cdot, ~\bar{}~, 0, 1)$ be an existentially closed De Morgan algebra. Then there is a definable (in the language of lattices) unary function $'$ on $D$ such that $(D;+,\cdot, ', ~\bar{}~, 0, 1)$ is a Boole-De Morgan algebra.
\end{corollary}

Now we prove that the theory of De Morgan (and Boole-De Morgan) algebras has the amalgamation property. This is actually well known (see \cite{cornish-fowler}) and follows from a more general theorem about amalgamation in varieties with certain properties. The same theorem can also be applied to Boole-De Morgan algebras. However, we give a proof below based on the amalgamation property of distributive lattices (see \cite{Gratzer-lattice-theory}).

\begin{proposition}
The theory of De Morgan algebras (Boole-De Morgan algebras) has the amalgamation property.
\end{proposition}
\begin{proof}
Let $D_0, D_1, D_2$ be De Morgan algebras with embeddings $f_i: D_0 \hookrightarrow D_i,~ i=1, 2$. By the amalgamation property of distributive lattices there is a distributive lattice $L$ such that $D_1$ and $D_2$ can be embedded into $L$ over $D_0$ as lattices. Now we embed each $D_i$ into $D_i \times D_i^{\op}$ as described above and notice that $D_1 \times D_1^{\op}$ and $D_2 \times D_2^{\op}$ can be embedded into $L \times L^{\op}$ over $D_0 \times D_0^{\op}$ (see the below diagram). These embeddings obviously respect the De Morgan negation and so $L \times L^{\op}$ is a De Morgan amalgam of $D_1$ and $D_2$ over $D_0$. 

A similar argument proves our claim for Boole-De Morgan algebras.

\[
  \begin{tikzcd}
  & D_1 \arrow{rr}{i_1} &  & D_1\times D_1^{\op} \arrow{dr}  & \\
  D_0 \arrow{rr}{i_0} \arrow{ur}{f_1} \arrow{dr}{f_2}  & &  D_0\times D_0^{\op} \arrow{ur} \arrow{dr} & & L\times L^{\op} \\
  & D_2 \arrow{rr}{i_2} & & D_2\times D_2^{\op} \arrow{ur}  & 
  \end{tikzcd}
\]
\end{proof}

\section{Existentially closed Boole-De Morgan algebras}

In order to characterise existentially closed Boole-De Morgan algebras we will prove that certain systems of equations and inequations always have a solution in such algebras. This will give rise to an \emph{Existential Closedness} axiom scheme. Then we will show that it is actually equivalent to being existentially closed.


\begin{definition}
Let $(A;+,\cdot, ', ~\bar{}~, 0, 1)$ be a Boole-De Morgan algebra. Define a map $ ^*:A\rightarrow A$ by $x^* = (\bar{x})' = \overline{(x')}$. It is evident that $ ^*$ is a lattice automorphism of order $2$. We extend the language of Boole-De Morgan algebras by adding a unary function symbol for $^*$.
\end{definition}

Assume $A$ is a finite Boole-De Morgan algebra with atoms $p_1,\ldots,p_n$. Then $^*$ permutes the atoms. Moreover the action of $^*$ on $A$ is completely determined by its action on $p_1,\ldots,p_n$. So it can be identified with a permutation $\sigma \in S_n$ (where $S_n$ is the symmetric group on the set $\{ 1, \ldots, n \}$), that is, $p_i^* = p_{\sigma(i)}$. Moreover, $\sigma$ must have order two, i.e. it is a product of disjoint two-cycles. We will write $\sigma_A$ for the permutation associated with $A$ and drop the subscript when no confusion can arise.

Now let $y:=(y_1,\ldots,y_n)$ be an $n$-tuple of variables. For a subset $I\subseteq \{ 1, \ldots, n\}$ consider the  formulae
\begin{align*}
    \psi^n_{I}(y,x) := & \bigwedge_{i\in I} y_i \cdot x\bar{x} = 0 \wedge \bigwedge_{i \notin I} y_i \cdot x\bar{x} \neq 0,\\
    \xi^n_{I}(y,x) := & \bigwedge_{i\in I} y_i \cdot xx^* = 0 \wedge \bigwedge_{i \notin I} y_i \cdot xx^* \neq 0,\\
    \chi^n_{I}(y,x) := & \bigwedge_{i\in I} y_i \cdot x'\bar{x} = 0 \wedge \bigwedge_{i \notin I} y_i \cdot x'\bar{x} \neq 0.
\end{align*}

For three sets $I_1,I_2,I_3 \subseteq \{ 1, \ldots, n\}$ we denote
$$\varphi^n_{(I_1,I_2,I_3)}(y,x) := \psi^n_{I_1}(y,x) \wedge \xi^n_{I_2}(y,x)\wedge \chi^n_{I_3}(y,x).$$


\begin{definition}
Let $\sigma \in S_n$ be a permutation. A triple $(I_1,I_2,I_3)$ of subsets of $\{ 1, \ldots, n\}$ is $\sigma$-\emph{consistent} if
\begin{itemize}
    \item[(i)]  $\sigma(I_2) = I_2,~ \sigma(I_3) = I_3$, and 
    \item[(ii)] $(I_1\cap I_2\cap I_3) \cap \sigma(I_1\cap I_2 \cap I_3) = \emptyset$.
\end{itemize}
\end{definition}

\begin{lemma}\label{consistent}
Let $A$ be a finite Boole-De Morgan algebra with atoms $p_1,\ldots,p_n,$ and $\sigma = \sigma_A$, that is, $p_i^* = p_{\sigma(i)}$ for all $i=1, \ldots,n$. If $\varphi^n_{(I_1,I_2,I_3)}(p,x)$ has a solution in an extension of $A$ (where $p:=(p_1,\ldots,p_n)$), then $(I_1,I_2,I_3)$ is $\sigma$-consistent.
\end{lemma}
\begin{proof}
Since $(xx^*)^* = xx^*$ and $(x'\bar{x})^* = x'\bar{x}$, the formulas $\xi^n_{I_2}(p,x)$ and $\chi^n_{I_3}(p,x)$ can have a realisation only if $I_2$ and $I_3$ are invariant under $\sigma$. Indeed, assume $\xi^n_{I_2}(p,u)$ holds for some $u$ in an extension of $A$. Then, using the fact that $\sigma^{-1} = \sigma$, we get $$i \in I_2 \mbox{ iff }p_i \cdot uu^* = 0 \mbox{ iff } p_{\sigma(i)} \cdot uu^* =0 \mbox{ iff } \sigma(i) \in I_2 \mbox{ iff } i\in \sigma(I_2).$$ Therefore, $I_2 = \sigma(I_2)$. Similarly, $I_3 = \sigma(I_3)$.

Further, assume there is an element $i\in (I_1\cap I_2\cap I_3) \cap \sigma(I_1\cap I_2 \cap I_3)$. Then $\sigma^{-1}(i) = \sigma(i) \in I_1\cap I_2\cap I_3$. If for some $u$ in an extension of $A$ the formula $\varphi^n_{(I_1,I_2,I_3)}(p,u)$ holds then $$p_i \cdot u\bar{u} = p^*_i \cdot u\bar{u} = p_i\cdot uu^* = p^*_i\cdot uu^* =p_i\cdot u'\bar{u}  = p^*_i\cdot u'\bar{u} = 0.$$
Therefore $p_i^* \cdot \bar{u} = p_i^* \cdot \bar{u} \cdot (u+u') = 0$ and so $p_i \cdot u' = 0$. On the other hand $p_i\cdot u = p_i \cdot u \cdot (\bar{u} + u^*) = 0$. Thus, $p_i = p_i \cdot (u+u') = 0$ which is a contradiction.
\end{proof}

\begin{definition}\label{defin-EC}
The theory $\ECBDA$ consists of the axioms of Boole-De Morgan algebras and  
the sentences
\begin{equation}
    \forall y_1, \ldots, y_n \left( \bigwedge_i y_i \neq 0 \wedge \sum_{i=1}^n y_i = 1 \wedge \bigwedge_{i \neq j} y_i \cdot y_j =0 \wedge \bigwedge_{i} y_i^* = y_{\sigma (i)} \rightarrow \exists x \varphi^n_{(I_1,I_2,I_3)}(y,x) \right) \tag{EC}
\end{equation} 
for each integer $n\geq 1$, each permutation $\sigma \in S_n$ with $\sigma^2 = \varepsilon$ (the identity permutation) and each $\sigma$-consistent triple $(I_1,I_2,I_3)$ of subsets of $\{ 1, \ldots, n\}$. 
\end{definition}

\begin{remark}
$\ECBDA$ stands for Existentially Closed Boole-De Morgan Algebras. We will prove shortly that ECBDA is indeed the theory of those algebras. The above axiom scheme is called \emph{Existential Closedness} (EC).
\end{remark}

\begin{remark}
It is easy to see that if $M$ is a Boole-De Morgan algebra satisfying EC for $n\leq 2$ then $M$ must be atomless, and hence infinite. Therefore ECBDA does not have finite models.
\end{remark}

\begin{theorem}\label{extension0}
Every Boole-De Morgan algebra can be extended to a model of $\ECBDA$.
\end{theorem}

First, we prove a lemma.

\begin{lemma}\label{lemma-I-J}
Assume $A\subseteq B$ are finite Boole-De Morgan algebras with atoms $p_1,\ldots,p_n$ and $q_1,\ldots, q_m$ respectively. Denote $p:=(p_1,\ldots,p_n),~ q:=(q_1,\ldots,q_m)$. Let $\sigma_A \in S_n,~ \sigma_B \in S_m$ be permutations with $p_i^* = p_{\sigma_A(i)},~ q_j^* = q_{\sigma_B(j)}$ for all $i, j$.  If for all $\sigma_B$-consistent  $I_1,I_2,I_3\subseteq \{ 1, \ldots, m\}$ the formula $\varphi^m_{(I_1,I_2,I_3)}(q,x)$ has a realisation in an extension of $B$ then for all $\sigma_A$-consistent $I_1,I_2,I_3 \subseteq \{ 1, \ldots, n\}$ the formula $\varphi^n_{(I_1,I_2,I_3)}(p,x)$ has a realisation in an extension of $A$.
\end{lemma}
\begin{proof} 
Assume $I_1,I_2,I_3 \subseteq \{ 1, \ldots, n\}$ are $\sigma_A$-consistent. Clearly each $p_i$ is the supremum of some $q_j$'s. For each $k=1,2,3$ denote $$J_k:= \{ j\in \{1,\ldots, m\}: q_j \leq p_i \mbox{ for some } i \in I_k \}.$$
We claim that $J_1,J_2,J_3$ are $\sigma_B$-consistent. Assume for contradiction that for some $l \in \{1,\ldots, m\}$ we have $l, \sigma_B(l) \in J_1\cap J_2 \cap J_3$ where $\sigma_B \in S_m$ with $q_j^* = q_{\sigma_B(j)}$. If for some $1\leq s,t \leq n$ we have $q_l\leq p_s,~ q_l^* \leq p_t$ then $q_l \leq p_sp_t^* = p_s p_{\sigma_A(t)}$. But if $s\neq \sigma_A(t)$ then $p_s p_{\sigma_A(t)} = 0$, hence $s= \sigma_A(t)$. Thus, for each $k$ there is $s_k \in I_k$ with $\sigma_A(s_k) \in I_k$ such that $q_l \leq p_{s_k}$. Then obviously $s_1=s_2=s_3=:s$. Thus, $s, \sigma_A(s) \in I_1\cap I_2\cap I_3$ which is a contradiction.

Now we show that a solution of $\varphi^m_{(J_1,J_2,J_3)}(q,x)$ is also a solution of $\varphi^n_{(I_1,I_2,I_3)}(p,x)$. Indeed, if $\varphi^m_{(J_1,J_2,J_3)}(q,u)$ holds for some $u$ in an extension of $B$ then by our definition of $J_k$ all equations in the system $\varphi^n_{(I_1,I_2,I_3)}(p,x)$ are satisfied at $x=u$. We need to prove that the inequations also hold at $u$. Pick $i\notin I_1$ and assume for contradiction that $p_i \cdot u\bar{u} = 0$. Then for every $q_j \leq p_i$ we must have $q_j \cdot u \bar{u} = 0$, hence $j \in J_1$. Therefore there is $i_j \in I_1$ such that $q_j\leq p_{i_j}$. Also, $i\neq i_j$ for $i \notin I_1$. Thus, $0 = p_i \cdot p_{i_j} \geq q_j$ which is a contradiction. The other inequations are dealt with similarly.
\end{proof}

\begin{proof}[Proof of Theorem \ref{extension0}]
Given a Boole-De Morgan algebra $A$ and a permutation $\sigma \in S_n$ with $\sigma^2= \varepsilon$, pick $p_1,\ldots,p_n\in A \setminus \{ 0 \}$ such that $p_i^* = p_{\sigma(i)}$, $\sum_{i=1}^n p_i = 1 $ and $p_i \cdot p_j =0$ whenever $i\neq j$. Take three $\sigma$-consistent sets $I_1,I_2,I_3$. We will show that in some extension of $A$ the formula $\varphi^n_{(I_1,I_2,I_3)}(p,x)$ has a solution. Then we can iterate this process and, taking the union of the obtained chain of structures, get an extension $A_1$ of $A_0:=A$ so that $A_1$ contains a realisation of $\varphi^n_{(I_1,I_2,I_3)}(p,x)$ for each appropriate choice of $n,~ \sigma, ~ I_1,~ I_2,~ I_3$ and $p_1,\ldots,p_n \in A$. Then we can construct $A_2, A_3, \ldots$ inductively where $A_{i+1}$ is the structure obtained from $A_i$ by the above procedure. Finally, the union $\bigcup_i A_i$ will be the desired extension of $A$ which is a model of ECBDA.

By the amalgamation property we may assume that $A$ is in fact the Boole-De Morgan algebra generated by $p_1,\ldots,p_n$. In particular, $A$ is finite and $p_1,\ldots,p_n$ are its atoms. 
Then by Proposition \ref{BDM-embedding-4} $A$ can be embedded into a direct power of $\mathbf{4}$, say $\mathbf{4}^m$. 
By the above lemma we may assume that $A=\mathbf{4}^m,~ n=2m$ and $p_i=(0,\ldots, a, \ldots, 0)$ (the $i$-th coordinate is $a$) for $1\leq i \leq m$ and  $p_i=(0,\ldots, b, \ldots, 0)$ (the $(i-m)$-th coordinate is $b$) for $m+1\leq i \leq 2m$ (these are all atoms of $\mathbf{4}^m$). Note also that in this case $\sigma(i) = m+i,~ 1\leq i \leq m$.

We show that there is an extension of $\mathbf{4}^m$ where $\varphi^n_{(I_1,I_2,I_3)}(p,x)$ has a solution. \\

\noindent \textbf{Case 1.} $m=1$. \\
In this case $p_1 = a,~ p_2 = b$ and $\sigma(1)=2$. We have only two possibilities for $I_2$ and $I_3$, either $\emptyset$ or $\{ 1, 2\}$. We embed $\mathbf{4}$ into $\mathbf{4}^k$ diagonally and show that in all cases a solution exists in the latter for some $k\leq 4$. 

\begin{itemize}
    \item If $I_1 = I_2 = \{ 1, 2\},~ I_3 = \emptyset$ then $x=0$ is a solution in $\mathbf{4}$.
    \item If $I_1 = \{ 1, 2\},~ I_2 = I_3 = \emptyset$ then $x=(1,0)$ is a solution in $\mathbf{4}^2$.
    
    \item If $I_1 = I_3 = \{ 1, 2\},~ I_2 = \emptyset$ then $x=1$ is a solution in $\mathbf{4}$.
    
    
    \item If $I_1 = \{ 1\},~ I_2 = I_3 = \emptyset$ then $x=(b,1,0)$ is a solution in $\mathbf{4}^3$.
    
    \item If $I_1 = \{ 1\},~ I_2 = I_3 = \{ 1, 2\}$ then $x=b$ is a solution in $\mathbf{4}$.
    
    \item If $I_1 = \{ 1\},~ I_2 = \{ 1, 2\},~ I_3 = \emptyset$ then $x=(b,0)$ is a solution in $\mathbf{4}^2$.
    
    \item If $I_1 = \{ 1\},~ I_2 =  \emptyset,~ I_3 = \{ 1, 2\}$ then $x=(b,1)$ is a solution in $\mathbf{4}^2$.
    
    \item If $I_1 = \{ 2\}$ then a solution can be found as in the previous four cases.
    
    \item If $I_1 =  I_2 = I_3 = \emptyset$ then $x=(a,b,0,1)$ is a solution in $\mathbf{4}^4$.
    
    \item If $I_1 = I_2 = \emptyset,~ I_3 = \{ 1, 2\}$ then $x=(a,b,1)$ is a solution in $\mathbf{4}^3$.
    
    \item If $I_1 =  I_3 = \emptyset,~ I_2 = \{ 1, 2\}$ then $x=(a,b,0)$ is a solution in $\mathbf{4}^3$.
    
    \item If $I_1 = \emptyset,~ I_2 = I_3 = \{ 1, 2\}$ then $x=(a,b)$ is a solution in $\mathbf{4}^2$\\
\end{itemize}

\noindent \textbf{Case 2.} $m>1$. \\
For $j=1,2,3$ and $1\leq i \leq m$ denote $I_j^i := I_j \cap \{ i, m+i \}$. Since $(I_1\cap I_2\cap I_3) \cap \sigma(I_1\cap I_2 \cap I_3) = \emptyset$, we cannot have $I_1^i = I_2^i = I_3^i = \{ i, m+i \}$. Hence for each $i$ there is an element $u_i$ in some extension $A_i$ of $\mathbf{4}$ such that $\varphi^2_{(I^i_1,I^i_2,I^i_3)}(a,b,u_i)$ holds. Then  $\varphi^n_{(I_1,I_2,I_3)}(p,x)$ is true of $(u_1,\ldots,u_m) \in A_1\times \ldots \times A_m$.
\end{proof}

\begin{notation}
For two Boole-De Morgan algebras $A\subseteq B$ and a subset $V:=\{ v_1,\ldots,v_n \}\subseteq  B$ the Boole-De Morgan subalgebra of $B$ generated by $A$ and $V$ is denoted by $A \langle V \rangle$ or $A\langle v_1,\ldots,v_n \rangle$.
\end{notation}

The following is an analogue of \cite[Lemma 7]{Schmid-ec-dl}.

\begin{lemma}\label{isomorphism}
Assume $A_0 \subseteq A \subseteq B$ are Boole-De Morgan algebras where $A\models \ECBDA$ and $A_0$ is finite. Then for every $v\in B$ there is $u\in A$ such that the map $f:A_0\cup \{ v \}  \rightarrow A_0 \cup \{ u \}$ which fixes $A_0$ pointwise and maps $v$ to $u$ extends to an isomorphism of $A_0\langle v \rangle$ and $A_0\langle u \rangle$.
\end{lemma}
\begin{proof}
Let $p_1, \ldots, p_n$ be the atoms of $A_0$ and denote $\sigma:= \sigma_{A_0}$. Define
\begin{align*}
    I_1  = & \{ i\in \{1,\ldots, n\} : p_i \cdot v\bar{v} = 0 \},\\
    I_2  = & \{ i\in \{1,\ldots, n\} : p_i \cdot vv^* = 0 \},\\
    I_3  = & \{ i\in \{1,\ldots, n\} : p_i \cdot v'\bar{v} = 0 \}.
\end{align*}

By Lemma \ref{consistent} the triple $(I_1,I_2,I_3)$ is $\sigma$-consistent. Let $u\in A$ be a realisation of the formula $\varphi^n_{(I_1,I_2,I_3)}(p,x)$.

The atoms of the Boole-De Morgan subalgebra of $B$ generated by $v$, which will be denoted by $C$, are in the set $\{ v\bar{v}, vv^*, v'\bar{v}, v'v^*, v, \bar{v}, v', v^*\}$. In particular, the last four elements can be expressed as joins of two of the first four elements. Hence, if any of the first four elements is non-zero then it is an atom. Also, notice that for any $p_i$ $$p_i \cdot v'v^* = 0 \Leftrightarrow p^*_i \cdot v\bar{v} = 0 \Leftrightarrow \sigma(i)\in I_1 \Leftrightarrow p^*_i \cdot u\bar{u} = 0 \Leftrightarrow p_i \cdot u'u^* = 0.$$

Thus, for any term function $g(x) \in \{ x\bar{x}, xx^*, x'\bar{x}, x'x^*, x, \bar{x}, x', x^*\}$ and any $p_i$ we have $$p_i \cdot g(v) =0 \mbox{ iff } p_i \cdot g(u) = 0.$$ This shows that the map $f$ fixing each $p_i$ and mapping $v$ to $u$ can be extended to a bijection between the atoms of $A_0 \langle v \rangle$ and $A_0 \langle u \rangle$ respecting the automorphism $^*$. Therefore, $f$ extends to an isomorphism of the underlying Boolean algebras of $A_0 \langle v \rangle$ and $A_0 \langle u \rangle$ which also respects $^*$ and hence the De Morgan negation. So it is actually an isomorphism of Boole-De Morgan algebras.
\end{proof}

\begin{theorem}\label{ec-ecbda}
A Boole-De Morgan algebra is existentially closed if and only if it is a model of $\ECBDA$.
\end{theorem}
\begin{proof}
By Theorem \ref{extension0} every existentially closed Boole-De Morgan algebra is a model of ECBDA. 

Now let $A \models \ECBDA$. Given a quantifier-free formula $\eta(x_1,\ldots,x_n)$ with parameters from a finite substructure $A_0 \subseteq A$ and a realisation $(v_1,\ldots,v_n)$ in an extension of $A$, we repeatedly apply Lemma \ref{isomorphism} and find $u_1,\ldots, u_n \in A$ such that $A_0\langle u_1,\ldots, u_n \rangle \cong A_0\langle v_1,\ldots, v_n \rangle$ with an isomorphism sending $u_i$ to $v_i$. Then $A \models \eta(u_1,\ldots, u_n)$.
\end{proof}

\begin{theorem}\label{BD-complete-qe}
$\ECBDA$ is the model completion of the theory of Boole-De Morgan algebras. It is complete and eliminates quantifiers.
\end{theorem}
\begin{proof}
Theorem \ref{ec-ecbda} shows that $\ECBDA$ is the model companion of the theory of Boole-De Morgan algebras. Since the latter is a universal theory and has the amalgamation property, the former is actually its model completion and admits quantifier elimination. Furthermore, the two-element Boole-De Morgan algebra $\mathbf{2}$ embeds into every Boole-De Morgan algebra, hence by quantifier elimination $\ECBDA$ is complete.
\end{proof}

\begin{example}
If ${B}$ is an atomless Boolean algebra then ${B} \times {B}^{\op}$ is an existentially closed Boole-De Morgan algebra and hence a model of $\ECBDA$.
\end{example}

\section{Existentially closed De Morgan algebras}

Now we translate the theory $\ECBDA$ to the language of De Morgan algebras replacing the Boolean negation by its definition in the language of lattices, that is,
$$z = x' \mbox{ iff } x+z = 1 \wedge x\cdot z = 0.$$

For $y:=(y_1,\ldots,y_n)$ and for a subset $I\subseteq \{ 1, \ldots, n\}$ consider the formulae 
\begin{align*}
    \tilde{\psi}^n_{I}(y,x) := & \bigwedge_{i\in I} y_i \cdot x\bar{x} = 0 \wedge \bigwedge_{i \notin I} y_i \cdot x\bar{x} \neq 0,\\
    \tilde{\xi}^n_{I}(y,x) := & \exists z \left( \bar{x} + z = 1 \wedge \bar{x}z =0 \wedge \bigwedge_{i\in I} y_i \cdot xz = 0 \wedge \bigwedge_{i \notin I} y_i \cdot xz \neq 0 \right),\\
    \tilde{\chi}^n_{I}(y,x) := & \exists z \left( x + z = 1 \wedge xz =0 \wedge \bigwedge_{i\in I} y_i \cdot z\bar{x} = 0 \wedge \bigwedge_{i \notin I} y_i \cdot z\bar{x} \neq 0 \right).
\end{align*}

Further, for $I_1,I_2,I_3 \subseteq \{ 1, \ldots, n\}$ set
$$\tilde{\varphi}^n_{(I_1,I_2,I_3)}(y,x) := \tilde{\psi}^n_{I_1}(y,x) \wedge \tilde{\xi}^n_{I_2}(y,x)\wedge \tilde{\chi}^n_{I_3}(y,x).$$

\begin{definition}
The theory ECDA consists of the axioms of complemented De Morgan algebras and  
the sentences
    $$\forall y_1,\ldots,y_n \left( \bigwedge_i y_i \neq 0 \wedge \sum_{i=1}^n y_i = 1 \wedge \bigwedge_{i \neq j} y_i y_j =0 \wedge \bigwedge_{i} \left(\bar{y}_i + y_{\sigma (i)} = 1 \wedge \bar{y}_i  y_{\sigma (i)} = 0 \right)  \rightarrow \exists x \tilde{\varphi}^n_{I}(y,x) \right)$$ for each integer $n\geq 1$, each permutation $\sigma \in S_n$ with $\sigma^2 = \varepsilon$ and each $\sigma$-consistent triple $I = (I_1,I_2,I_3)$ of subsets of $\{ 1, \ldots, n\}$.
\end{definition}

\begin{theorem}
$\ECDA$ is the model completion of the theory of De Morgan algebras. It is complete and eliminates quantifiers.
\end{theorem}
\begin{proof}
This follows from the results of the previous two sections.
\end{proof}

\begin{example}
If ${L}$ is a bounded atomless complemented distributive lattice then ${L} \times {L}^{\op}$ is a model of $\ECDA$.
\end{example}

\section{Model theoretic properties}

\subsection{$\aleph_0$-categoricity} The theory of atomless Boolean algebras is $\aleph_0$-\emph{categorical}, that is, it has a unique countable model up to isomorphism. We show now that $\ECBDA$ and $\ECDA$ have the same property.

\begin{theorem}\label{cat}
$\ECBDA$ and $\ECDA$ are $\aleph_0$-categorical.
\end{theorem}

The proof uses the Ryll-Nardzewski theorem which is recalled below for the convenience of the reader. 

\begin{theorem}[Ryll-Nardzewski, {\cite[Theorem 4.4.1]{Mar}}]
A theory $T$ is $\aleph_0$-categorical if and only if for every positive integer $n$ there are only finitely many formulas with $n$ free variables modulo $T$, that is, every formula with $n$ free variables is equivalent to one from a fixed finite set of formulas.
\end{theorem}

\begin{proof}[Proof of Theorem \ref{cat}]
It suffices to prove this for $\ECBDA$. The proof is based on local finiteness of Boole-De Morgan algebras (Corollary \ref{local-finite}). For each positive integer $n$ the free $n$-generated Boole-De Morgan algebra is finite (cf. \cite[Section 4]{Mov-Asl-Boole-de}). So there are finitely many terms $t_i(x_1,\ldots,x_n),~ i\in I$ such that for every term $t(x_1,\ldots,x_n)$ there is an $i\in I$ with
$$\ECBDA \models \forall x_1,\ldots,x_n (t(x_1,\ldots,x_n) = t_i(x_1,\ldots,x_n)).$$

By quantifier elimination, every formula is a Boolean combination of formulas of the form $t_1= t_2$ or $t_1 \neq t_2$ where $t_1$ and $t_2$ are terms. So the above observation implies that for every $n$ there are finitely many formulas of $n$ variables modulo $\ECBDA$. Now $\aleph_0$-categoricity follows from the Ryll-Nardzewski theorem.
\end{proof}

The above argument also shows that the theories of atomless Boolean algebras and atomless complemented distributive lattices are $\aleph_0$-categorical. 

\begin{remark}
The unique countable model of ECDA (ECBDA) is the Fra\"iss\'e limit of finite De Morgan (respectively Boole-De Morgan) algebras, and is homogeneous. Homogeneity, as well as $\aleph_0$-categoricity, can be deduced directly from Lemma \ref{isomorphism} by a back-and-forth argument.
\end{remark}

\begin{remark}
It is clear that models of $\ECDA$ and $\ECBDA$ are atomless. This means that we have infinite linearly ordered sets in all models of these theories. Hence, for each uncountable cardinal $\kappa$ these theories have $2^{\kappa}$ non-isomorphic models of cardinality $\kappa$ (see \cite[Chapter 5, Theorem 5.3.2]{Mar}).
\end{remark}

\subsection{Algebraic and definable closures}

In this section we describe algebraic and definable closures in $\ECBDA$ and $\ECDA$. We begin by recalling the necessary model theoretic notions and fixing some notation. As before, the reader is referred to \cite{Mar,TZ} for details.

\begin{definition}
Let $M$ be a structure and $A\subseteq M$ be a subset (possibly empty). The \emph{algebraic (definable) closure} of $A$ in $M$, denoted by $\acl(A)$ (respectively $\dcl(A)$), is the union of all definable finite sets (respectively, singletons) with parameters from $A$.
\end{definition}

For example, in an algebraically closed field the model theoretic definable closure of a set $A$ is the subfield generated by $A$, and the model theoretic algebraic closure of $A$ coincides with the field theoretic algebraic closure of that subfield. The definable and algebraic closures in atomless Boolean algebras and atomless complemented bounded distributive lattices coincide with the Boolean subalgebra generated by the parameter set.

\begin{definition}
For a structure $M$, a subset $A \subseteq M$ and a tuple $\bar{b}\in M^n$ the \emph{type of} $\bar{b}$ over $A$ (in $M$), denoted $\tp(\bar{b}/A)$, is the set of all formulas with parameters from $A$ that are true of $\bar{b}$ in $M$. A \emph{realisation} of a type is a tuple which satisfies all formulas of the type. A type is \emph{isolated} if there is a formula in the type which implies all other formulas of the type. A type is \emph{algebraic} if it has only finitely many \emph{realisations}.
\end{definition}

\begin{notation}
For a positive integer $n$ the set $\{ 1,\ldots, n\}$ is denoted by $[n]$.
\end{notation}

\begin{theorem}\label{acl-dcl}
If $D$ is a model of $\ECBDA$ or $\ECDA$ then the definable and algebraic closures of an arbitrary set $A\subseteq D$ are both equal to the Boole-De Morgan subalgebra of $D$ generated by $A$.
\end{theorem}

We will prove this for $\ECBDA$ only. Assume $D \models \ECBDA$ and $A\subseteq D$ is a finite Boole-De Morgan subalgebra of $D$ with atoms $p_1, \ldots, p_n$. Pick an element $v \in D$ and denote 
\begin{align*}
    I_1  = & \{ i\in [n] : p_i \cdot v\bar{v} = 0 \},\\
    I_2  = & \{ i\in [n] : p_i \cdot vv^* = 0 \},\\
    I_3  = & \{ i\in [n] : p_i \cdot v'\bar{v} = 0 \}.
\end{align*}

\begin{lemma}
The formula $\varphi^n_{(I_1,I_2,I_3)}(p,x)$ isolates the type $\tp(v/A)$, that is, for any formula $\psi(x)$ from that type (with parameters from $A$) we have $D \models \forall x (\varphi^n_{(I_1,I_2,I_3)}(p,x) \rightarrow \psi(x))$.
\end{lemma}
\begin{proof}
This follows from quantifier elimination and the proof of Lemma \ref{isomorphism}.
\end{proof}

\begin{lemma}
Assume for some set $I\subseteq  [n]$ we have $v=\sum_{i\in I} p_i$. Then 
\begin{gather*}
    I_1 = ([n]\setminus I) \cup \sigma_A(I),\\
    I_2 = [n]\setminus (I \cap \sigma_A(I)),\\
    I_3 = I \cup \sigma_A(I).
\end{gather*}
Conversely, if the above equalities hold for some $I$ then $v=\sum_{i\in I} p_i$ and it is the only realisation of $\varphi^n_{(I_1,I_2,I_3)}(p,x)$ in an extension of $D$.
\end{lemma}
\begin{proof}
Assume $v=\sum_{i\in I} p_i$. Let us prove the first equality. The other two equalities are proven similarly. To this end we notice that
$$i\in I_1 \mbox{ iff } p_i \cdot v\bar{v} = 0 \mbox{ iff } (p_i \cdot v = 0 \mbox{ or } p_i \cdot \bar{v} = 0) \mbox{ iff } (i\notin I \mbox{ or } i\in \sigma_A(I)).$$
Now if the equalities hold for some $i$ then the element $u:=\sum_{i\in I} p_i$ realises the formula $\varphi^n_{(I_1,I_2,I_3)}(p,x)$ and hence the type $\tp(v/A)$. But since $u\in A$, it is the unique realisation of $\tp(v/A)$. In particular, $v=u$.
\end{proof}

\begin{definition}
We say $(I_1,I_2,I_3)$ is $A$-trivial (or just trivial) if the above equalities hold for some $I\subseteq [n]$.
\end{definition}

Triviality just means that the formula $\varphi^n_{(I_1,I_2,I_3)}(p,x)$ is equivalent to $x=u$ for some $u\in A$.

\begin{lemma}\label{extension}
Let $A\subseteq B \subseteq D$ where $B$ is a finite Boole-De Morgan extension of $A$ with atoms $q_1,\ldots,q_m$. For $k=1,2,3$ denote $$J_k:= \{ j\in[m]: q_j \leq p_i \mbox{ for some } i \in I_k \}.$$
If $(J_1,J_2,J_3)$ is $B$-trivial then $(I_1,I_2,I_3)$ is $A$-trivial.
\end{lemma}
\begin{proof}
For each $1\leq i \leq n$ denote $Q_i:= \{ j\in [m]: q_j \leq p_i\}$. Then $Q_i$'s are pairwise disjoint, $J_k = \bigcup_{i\in I_k}Q_i$ and $$\sigma_A(i_1) = i_2 \mbox{ iff } \sigma_B(Q_{i_1}) = Q_{i_2}. $$

Now assume $J\subseteq [m]$ witnesses the triviality of $(J_1,J_2,J_3)$. We claim that $$I:= \{ i \in [n]: Q_i \subseteq J \} $$ witnesses the triviality of $(I_1,I_2,I_3)$.

We show first that $J$ is a union of some $Q_i$'s. Denote $\tilde{J}:= J \setminus \bigcup_{i\in I} Q_i$. We will show that $\tilde{J} = \emptyset$. Since $J_1 = ([m]\setminus J) \cup \sigma_B(J)$ is a union of some $Q_i$'s, $\sigma_B(J) \supseteq \tilde{J}$. On the other hand $J_2 = [m] \setminus (J \cap \sigma_B(J))$ is also a union of some $Q_i$'s, therefore so is $J \cap \sigma_B(J)$. The latter is equal to $$ \left( \bigcup_{i\in I} Q_i \cap \sigma_B(J) \right) \cup (\tilde{J}\cap \sigma_B(J)) = \left( \bigcup_{i\in I} Q_i \cap \sigma_B(J) \right) \cup \tilde{J}.$$ This implies $\tilde{J}=\emptyset$ for $\bigcup_{i\in I} Q_i$ is disjoint from $\tilde{J}$.

Now we show that $(I_1,I_2,I_3)$ is trivial. Let $i \in [n]$.
\begin{itemize}
    \item We have $i\in I_1$ iff $Q_i \subseteq J_1 = ([m]\setminus J) \cup \sigma_B(J)$ iff $Q_i \subseteq [m] \setminus J$ or $Q_i \subseteq \sigma_B(J)$. The former is equivalent to $i \in [n]\setminus I$, while the latter is the case iff $Q_{\sigma_A(i)} \subseteq J$ iff $\sigma_A(i) \in I$ iff $i \in \sigma_A(I)$. Thus $I_1 = ([n]\setminus I) \cup \sigma_A(I)$.
    
    \item  We have $i\in I_2$ iff $Q_i \subseteq J_2 = [m]\setminus (J \cap \sigma_B(J))$ iff $Q_i \cap J \cap \sigma_B(J) = \emptyset$ iff $Q_i \cap J = \emptyset$ or $Q_i \cap \sigma_B(J) = \emptyset$. As above, this is equivalent to $i \notin I$ or $i \notin \sigma_A(I)$, hence $I_1 = [n]\setminus (I \cap \sigma_A(I))$.
    
      \item We have $i\in I_3$ iff $Q_i \subseteq J_3 = J \cup \sigma_B(J)$ iff $Q_i \subseteq J$ or $Q_i \subseteq \sigma_B(J)$. This happens iff $i\in I$ or $i\in \sigma_A(I)$ and so $I_3 = I \cup \sigma_A(I)$.
\end{itemize}

\end{proof}

Now we are ready to prove the main theorem. 

\begin{proof}[Proof of Theorem \ref{acl-dcl}]

Assume $A \subseteq D$ is a finite Boole-De Morgan subalgebra with atoms $p_1, \ldots, p_n$. Pick an element $v \in D\setminus A$ and let $I_1, I_2, I_3$ be as above. Then $(I_1,I_2,I_3)$ is non-trivial. Let $B := A\langle v \rangle$ be the Boole-De Morgan subalgebra of $D$ generated by $A\cup \{ v \}$ and let $q_1,\ldots,q_m$ be the atoms of $B$. Define $J_1,J_2,J_3$ as in Lemma \ref{extension}. Then $(J_1,J_2,J_3)$ is non-trivial. By the EC axiom scheme there is an element $w\in D$ such that $\varphi^m_{(J_1,J_2,J_3)}(q,w)$. By the proof of Lemma \ref{lemma-I-J}, $w$ also satisfies the formula $\varphi^n_{(I_1,I_2,I_3)}(p,x)$, hence it realises the type $\tp(v/A)$. On the other hand $w\notin A$ for $(J_1,J_2,J_3)$ is non-trivial. Thus, $\tp(v/A)$ has a realisation $w\neq v$. Repeating this procedure, we will find infinitely many realisations of $\tp(v/A)$, hence $v \notin \acl(A)$.
\end{proof}

\subsection*{Acknowledgements} I am grateful to James Raftery for referring me to \cite{Clark-Davey} for a discussion of existentially closed algebras from the point of view of natural dualities. I would also like to thank the referee for numerous useful remarks which helped to improve the presentation of the paper.

\addcontentsline {toc} {section} {Bibliography}
\bibliographystyle {alpha}
\bibliography {ref}

\begin{thebibliography}{ABD92}

\bibitem[AB75]{Ross-Belnap}
Alan~Ross Anderson and Nuel~D. Belnap, Jr.
\newblock {\em Entailment. {T}he logic of relevance and necessity. {V}ol. {I}}.
\newblock Princeton University Press, Princeton, N. J.-London, 1975.

\bibitem[ABD92]{Ross-belnap-dunn}
Alan~Ross Anderson, Nuel~D. Belnap, Jr., and J.~Michael Dunn.
\newblock {\em Entailment. {T}he logic of relevance and necessity. {V}ol.
  {II}}.
\newblock Princeton University Press, Princeton, NJ, 1992.

\bibitem[BC14]{twist-product}
Manuela Busaniche and Roberto Cignoli.
\newblock The subvariety of commutative residuated lattices represented by
  twist-products.
\newblock {\em Algebra Universalis}, 71(1):5--22, 2014.

\bibitem[BD74]{balbes-dwinger}
Raymond Balbes and Philip Dwinger.
\newblock {\em Distributive lattices}.
\newblock University of Missouri Press, Columbia, Mo., 1974.

\bibitem[BEI03]{Brzozowski-hazard}
Janusz Brzozowski, Zolt\'{a}n \'{E}sik, and Yaacov Iland.
\newblock Algebras for hazard detection.
\newblock In {\em Beyond two: theory and applications of multiple-valued
  logic}, volume 114 of {\em Stud. Fuzziness Soft Comput.}, pages 3--24.
  Physica, Heidelberg, 2003.

\bibitem[Bir67]{birkhoff}
Garrett Birkhoff.
\newblock {\em Lattice theory}.
\newblock Third edition. American Mathematical Society Colloquium Publications,
  Vol. XXV. American Mathematical Society, Providence, R.I., 1967.

\bibitem[CD98]{Clark-Davey}
David~M. Clark and Brian~A. Davey.
\newblock {\em Natural dualities for the working algebraist}, volume~57 of {\em
  Cambridge Studies in Advanced Mathematics}.
\newblock Cambridge University Press, Cambridge, 1998.

\bibitem[CF77]{cornish-fowler}
William~H. Cornish and Peter~R. Fowler.
\newblock Coproducts of {D}e {M}organ algebras.
\newblock {\em Bull. Austral. Math. Soc.}, 16(1):1--13, 1977.

\bibitem[Gr{\"{a}}11]{Gratzer-lattice-theory}
George Gr{\"{a}}tzer.
\newblock {\em Lattice Theory: Foundation}.
\newblock Birkh\"{a}user/Springer Basel AG, Basel, 2011.

\bibitem[H{\'{a}}j98]{hajek-fuzzy}
Petr H{\'{a}}jek.
\newblock {\em Metamathematics of fuzzy logic}, volume~4 of {\em Trends in
  Logic---Studia Logica Library}.
\newblock Kluwer Academic Publishers, Dordrecht, 1998.

\bibitem[Kal58]{kalman}
J.~A. Kalman.
\newblock Lattices with involution.
\newblock {\em Trans. Amer. Math. Soc.}, 87:485--491, 1958.

\bibitem[MA14]{Mov-Asl-Boole-de}
Yuri Movsisyan and Vahagn Aslanyan.
\newblock Boole-{D}e {M}organ algebras and quasi-{D}e {M}organ functions.
\newblock {\em Communications in Algebra}, 42(11):4757--4777, 2014.

\bibitem[Mar02]{Mar}
David Marker.
\newblock {\em Model Theory: An Introduction}.
\newblock Springer, 2002.

\bibitem[Moi35]{moisil}
Grigore Moisil.
\newblock Recherches sur l'alg\'ebre de la logique.
\newblock {\em Annales scientifiques de l'Universit\'e de Jassy}, 22:1--117,
  1935.

\bibitem[Poi00]{Poizat}
Bruno Poizat.
\newblock {\em A Course in Model Theory}.
\newblock Springer, 2000.

\bibitem[Sch79]{Schmid-ec-dl}
J\"{u}rg Schmid.
\newblock Algebraically and existentially closed distributive lattices.
\newblock {\em Z. Math. Logik Grundlag. Math.}, 25(6):525--530, 1979.

\bibitem[TZ12]{TZ}
Katrin Tent and Martin Ziegler.
\newblock {\em A Course in Model Theory}.
\newblock Cambridge University Press, 2012.

\end{thebibliography}

\end{document}